\documentclass[11pt,leqno,fleqn]{article}
\usepackage{latexsym,amsmath,amssymb}
\usepackage[dvips]{epsfig}
\usepackage{epsfig}

\newcommand{\AS}{\normalsize}

\newcommand{\lba}{\left(\small\begin{array}{c}}
\newcommand{\ear}{\end{array}\right)}
\newcommand{\ee}{\end{equation}}
\newcommand{\bea}{\begin{eqnarray}}
\newcommand{\eea}{\end{eqnarray}}
\newcommand{\bean}{\begin{eqnarray*}}
\newcommand{\eean}{\end{eqnarray*}}

\newcommand{\ecke}{\;_-\!\rule{0.2mm}{0.2cm}\;\;}

\newcommand{\Cal}{\mathcal}
\newcommand{\SO}{SO}

\newcommand{\C}{\Bbb{C}}
\newcommand{\R}{\Bbb{R}}

\newtheorem{theorem}{Theorem}[section]
\newtheorem{lemma}[theorem]{Lemma}
\newtheorem{prop}{Proposition}[section]
\newtheorem{co}[theorem]{Corollary}

\newtheorem{definition}[theorem]{Definition}
\newtheorem{example}[theorem]{Example}

\newtheorem{remark}[theorem]{Remark}

\begin{document}
\title{Imaginary Killing Spinors in Lorentzian Geometry}
\author{Felipe Leitner
\footnote{The author was supported by EPSRC Grant GR/R62694.}}
\date{January 2003}
\maketitle
\begin{sloppypar}
\begin{abstract}
We study the geometric structure of Lorentzian spin
manifolds, which admit imaginary Killing spinors. The discussion
is based on the cone construction and a normal form classification
of skew-adjoint operators in signature $(2,n-2)$. Derived
geometries include Brinkmann spaces, Lorentzian Einstein-Sasaki spaces
and certain warped product structures. Exceptional cases with decomposable 
holonomy of the cone are possible.   
\end{abstract}
\section{Introduction}
A classical object of interest in differential geometry are Killing vector
fields. These are by definition infinitesimal isometries, which means
that the flow of such a vector field preserves the metric. A spinorial
analog 
are the so-called Killing spinor fields, which occur on spin
manifolds and are defined as solutions of the field equation
$\nabla^S_X \varphi =\lambda X \cdot\varphi$ for all vector fields $X$
and some fixed $\lambda\in\C$,
where $\nabla^S$ denotes the spinor
derivative and $\cdot$ the Clifford
multiplication.\\ 

In Riemannian geometry, it was proved in \cite{Fri80} that real Killing
spinors
realize the lower bound of the eigenvalue estimation for the
Dirac equation on compact spaces with positive scalar curvature.
In the sequel Riemannian spaces admitting Killing spinors were
intensively studied (cf. \cite{DNP86}, \cite{BFGK91}) and a
complete
geometric
description of such spaces was established. For the case of 
imaginary Killing spinors this was done by H. Baum in \cite{Bau89}
and then for real Killing spinors by Ch. B{\"a}r using
the cone construction and the holonomy classification
of Riemannian spaces with parallel spinors (cf. \cite{Bar93}). Both
results
characterize Riemannian spaces with Killing spinors by the Einstein
condition and the existence of certain differential forms, which can be 
understood as generalized Killing vectors. Real Killing spinors in
Lorenztian geometry were first studied in \cite{Boh98}.\\

In this paper we will tread the Killing spinor equation to an
imaginary Killing number $\lambda$ 
on a pseudo-Riemannian space with Lorentzian signature. 
As technical tool we will use again the cone construction for the
investigation. 
Contrary to the Riemannian case, a holonomy description of the cone can
not be used, since there
is no classification of indecomposable holonomy groups
for pseudo-Riemannian manifolds. 
Moreover, the geodesical completness of a Lorenztian manifold
does not imply that a cone with decomposable holonomy is flat. 
Instead, our geometrical description is mainly based
on a normal form classification of skew-adjoint operators
in signature $(2,n-2)$, which is more rich then in the Euclidean
case (cf. \cite{Bou00}). The derived Lorenztian geometries are then
described
by the causal properties of the corresponding Dirac current and 
the existence of parallel spinors or certain Killing forms.
Thereby, we will use the knowledge of structure results for
Lorentzian manifolds admitting conformal gradient fields (cf.
\cite{KR97}) and twistor spinors with lightlike Dirac current (cf.
\cite{BL02}). 
Examples of geometries that occur are the Brinkmann spaces with parallel
spinors and the Lorenztian Einstein-Sasaki manifolds.\\

The order of the paper is as follows. In the next section
we introduce the basic notations and definitions appropriate
for the study of Killing spinors and state basic
curvature conditions for their existence (cf. Proposition
\ref{p21} and \ref{p22}). In section 3 we
recall
the cone construction over a Lorenztian base manifold and the
correspondence of 'Killing objects'
on the base and parallel objects on its cone (Theorem \ref{t31}). 
We
present
the normal form classification of skew-adjoint operators
in signature $(2,n-2)$ due to the work of Ch. Boubel 
in section 4. It turns out that there
are exactly four generic types of normal forms for skew-adjoint 
operators coming from a spinor (cf. Corollary \ref{c47}).
The cone of a Lorentzian manifold admitting imaginary Killing 
spinors is furnished with at least one parallel $2$-form, which
corresponds to one of the generic types (Proposition \ref{p51}). 
According to this type of a parallel $2$-form on the cone
we undertake in three of the four generic cases a discussion of the
geometry of Lorenztian manifolds
with imaginary Killing spinors in the last section. 
The results of the discussion are summerized in Theorem \ref{t53}.

\section{Basic facts on Killing spinors}  %

In this section we recall the definition of Killing spinors on a spin 
manifold and fix some notations. A basic integrability condition for
Killing spinors is stated. For more details we refer to \cite{BFGK91}. 
Moreover, we will come across special Killing forms as they were
introduced in \cite{Sem02}.\\
 
Let $(M^{n,k},g)$ be a semi-Riemannian spin manifold of dimension $n\geq
3$ and signature $(k,n-k)$ ($k$ is the number of timelike vectors in
an orthonormal basis at a point). We denote by $S$ the complex spinor
bundle and by $\cdot$ the Clifford multiplication on spinors. The Dirac 
operator $D:\Gamma(S)\rightarrow\Gamma(S)$ acting on smooth spinor fields
is defined as superposition of spinor derivative $\nabla^S$ and Clifford
multiplication. 
A spinor field $\varphi\in\Gamma(S)$ is called {\it{Killing spinor}} to
the
{\it{Killing number}} $\lambda\in\C$ if it satisfies the equation
\[
\nabla^S_X\varphi=\lambda X\cdot\varphi\qquad 
\mathrm{for\ all\ vector\ fields}\ \ X.
\]
It follows immediately from this definition that a Killing spinor
$\varphi$
is
an eigenspinor of the Dirac operator $D$ to the eigenvalue $-n\lambda$ 
and $\varphi$ is obviously a parallel spinor field with respect to the
modified 
spinor derivative $\tilde{\nabla}_\lambda$ defined by
\[ 
\tilde{\nabla}_\lambda:=\nabla^S-\lambda id_{TM}.
\] 
In particular, this implies that a Killing spinor $\varphi$ admits no
zeros. It holds the following basic integrability condition.

\begin{prop} \label{p21} (\cite{BFGK91}) Let $\varphi\in\Gamma(S)$ be
a Killing
spinor to the Killing number $\lambda\in\C$. 
\begin{enumerate}
\item It is $\Cal{W}(\eta)\cdot\varphi=0$ for any
$2$-form $\eta$, where $\Cal{W}$ denotes the Weyl tensor.
\item 
$(Ric(X)-4\lambda^2(n-1)X)\cdot\varphi=0$, i.e. the image 
of the map $Ric-4\lambda^2(n-1)id_{TM}$ is totally lightlike or trivial.
\item
The scalar curvature is constant and given by $scal=4n(n-1)\lambda^2$.
The Killing number $\lambda$ is real or purely imaginary.
\end{enumerate}
\end{prop}

If the Killing number $\lambda$ is zero ($scal=0$), $\varphi$ is a
parallel spinor, in case that $\lambda$ is real and non-zero ($scal >0$),
$\varphi$ is called real Killing spinor, and in case that $\lambda$ is
purely imaginary ($scal < 0$), $\varphi$ is called imaginary Killing
spinor. We will treat in this paper the Killing spinor equation with 
imaginary 
Killing number on a space of Lorentzian signature 
$(-+\ldots+)$.\\

Let $(M^{n,1},g)$ 
be a connected, oriented and time-oriented Lorentzian spin manifold. 
There exists an indefinite non-degenerate inner product
$\langle\cdot,\cdot\rangle$ on the spinor
bundle $S$
such that
\begin{eqnarray*}  
\langle X \cdot \varphi  , \psi \rangle & =& \langle \varphi , X \cdot
\psi
\rangle \quad \quad \mbox{and}\\
X(\langle \varphi , \psi \rangle )&= &\langle \nabla^S_X \varphi , \psi
\rangle + \langle \varphi , \nabla^S_X \psi \rangle
\end{eqnarray*}
for all vector fields $X$ and all spinor fields $\varphi,\psi$. Each
spinor field $\varphi\in\Gamma(S)$ defines a vector field $V_\varphi$
on $M$, the so-called {\it Dirac current}, by the relation
$g(V_\varphi,X):=-\langle X\cdot\varphi,\varphi\rangle$ for all
vector fields $X$.
The Dirac current satisfies the following pointwise properties.

\begin{lemma} \label{l21} (\cite{Lei01}) Let $(M^{n,1},g)$ be a Lorentzian
spin manifold and let $\varphi(p)\neq 0$ be a spinor in a point $p\in
M^{n,1}$. Then
\begin{enumerate}
\item $V_\varphi(p)\neq 0\ $ and\
$V_\varphi(p)$ is causal (i.e. $g_p(V_\varphi,V_\varphi)\leq 0$).
\item
If $X\!\cdot\varphi(p)=\rho\varphi(p)$ for some $0\neq\! X\!\in\!
T_pM$
and
$\rho\in\!\R$ then the vector $X$ is parallel to $V_\varphi(p)$.
\end{enumerate}
\end{lemma}

The lemma makes clear that the Dirac current to a Killing
spinor on a Lorentzian manifold is everywhere causal. Moreover,
it is now possible to prove a stronger 
curvature condition for the existence of Killing spinors.

\begin{prop} \label{p22} Let $(M^{n,1},g)$ be a Lorentzian spin manifold 
admitting a Killing spinor $\varphi$, whose Dirac current $V_{\varphi}$
is timelike. Then $(M,g)$ is an Einstein space.
\end{prop}

{\sc Proof}: Let us assume that $M^{n,1}$ is a non-Einstein space.
Then there is an open set $U$ in $M$,
where $H:=Ric(X)-\frac{scal}{n}X\neq 0$ is lightlike for some vector field
$X$.
The Clifford product $H\cdot\varphi$ vanishes and, by Lemma \ref{l21},
this
implies that $H$ and 
$V_\varphi$ are parallel, which is a
contradiction to the assumption.\hfill$\Box$\\

Especially, for imaginary Killing spinors it holds

\begin{prop} \label{p23} (\cite{BFGK91}) Let $\varphi$ be an imaginary
Killing spinor 
on a Lorenztian spin manifold $(M^{n,1},g)$. Then
the length $\langle\varphi,\varphi\rangle$ is constant on $M^{n,1}$
and if $\langle\varphi,\varphi\rangle\not\equiv 0$ the space
$M^{n,1}$ is Einstein.
\end{prop}

{\sc Proof}: It is $X(\langle\varphi,\varphi\rangle)=
\langle\lambda X\cdot\varphi,\varphi\rangle+\langle\varphi,
\lambda X\cdot\varphi\rangle=0$ and
with Proposition \ref{p21} we calculate 
\begin{eqnarray*}\frac{1}{4(n-1)}Ric(X,Y)\langle\varphi,\varphi\rangle
&=&\frac{-1}{4(n-1)}Re\langle Ric(X)\cdot\varphi,Y\cdot\varphi\rangle
=-Re\langle\lambda^2 X\cdot\varphi,Y\cdot\varphi\rangle\\[2mm]&=&
\lambda^2
g(X,Y)\langle\varphi,\varphi\rangle\end{eqnarray*}
for all vector fields $X$ and $Y$, which shows that
$Ric(X)=\frac{scal}{n}X$ in case that  $\langle\varphi,\varphi\rangle\neq
0$.\hfill$\Box$\\

Proposition \ref{p22} and \ref{p23} imply that an imaginary
Killing spinor $\varphi$
on a Lorentzian non-Einstein space $M^{n,1}$ must have
vanishing length 
$\langle\varphi,\varphi\rangle\equiv 0$ and the Dirac current $V_\varphi$
to $\varphi$ must be lightlike on an open subset of $M^{n,1}$. We will
see later that in this case $V_\varphi$ is even lightlike everywhere on 
$M^{n,1}$. Moreover, the Dirac current satisfies

\begin{prop} \label{p24} Let $\varphi$ be an imaginary Killing spinor 
on a Lorentzian spin manifold $(M^{n,1},g)$. 
The Dirac current $V_\varphi$ is a Killing vector field, which 
in addition satisfies 
$\nabla_XdV_{\varphi}^{\flat}=-4\lambda^2X^{\flat}\wedge
V_{\varphi}^{\flat}$.
\end{prop}

{\sc Proof}: It holds 
\[ g(\nabla_{e_i}V_\varphi,e_j)=
-\langle\lambda e_je_i\cdot\varphi,\varphi\rangle
-\langle e_j\varphi,\lambda e_i\cdot\varphi\rangle
=-g(\nabla_{e_j}V_\varphi,e_i)\]
for all $i,j\in\{1,\ldots ,n\}$, where $(e_1,\ldots,e_n)$ is an arbitrary 
orthonormal basis on $M^{n,1}$. This proves that $V_\varphi$ is a 
Killing 
vector
field. Moreover, 
\[\begin{array}{l}dV_\varphi^b=
4\lambda\sum_{i<j}\varepsilon_i\varepsilon_j\langle 
e_ie_j\cdot\varphi,\varphi\rangle
e_i^\flat\wedge e_j^\flat\qquad\qquad \mbox{and} \\[3mm]
\nabla_XdV_\varphi^\flat=4\lambda^2\sum_{i<j}\varepsilon_i
\varepsilon_j(\langle e_ie_jX\cdot\varphi,\varphi\rangle-
\langle Xe_ie_j\cdot\varphi,\varphi\rangle)e_i^\flat\wedge
e_j^\flat=-4\lambda^2
X^\flat\wedge V^\flat_\varphi\ .\end{array}\]
\hfill$\Box$\\

In general, a $p$-form $\alpha^p$, which solves the equation 
\[\nabla_Xd\alpha^p=cX^{\flat}\wedge\alpha^p\qquad\quad\mbox{for\ all\
vectors}\ X\]
and some fixed
$c\in\R$, is called a
{\it special Killing $p$-form} (cf. \cite{Sem02}). Proposition \ref{p24} 
states that
the dual of the Dirac current to an imaginary Killing spinor 
is a special Killing $1$-form. 
Killing spinors also produce special Killing forms 
of other degree then $1$. For this, we observe that one
constructs a $p$-form $\alpha^p_\varphi$ to a spinor $\varphi$ by the 
rule \[g(\alpha^p_\varphi,X^p):=
-i^{\frac{p(p-1)}{2}}\langle X^p\cdot\varphi,\varphi\rangle\qquad
\mbox{for}\ \mbox{all}\ p\mbox{-forms}\ X^p\]
and, in fact, if $p$ is odd and the Killing number $\lambda$ of a 
Killing spinor $\varphi$ is 
imaginary (or if  $p$ is even 
and $\lambda$ is real) then the
associated $p$-form
$\alpha^p_\varphi$ to $\varphi$ is special Killing.
      
\section{The cone $\hat{M}$}   %

We defined in the last section Killing spinors and special
Killing $p$-forms on a Lorenztian manifold. In this section we will 
interpret these as parallel objects on the {\it cone manifold}. 
The cone construction was originally applied in order to 
describe Riemannian geometries admitting real Killing spinors (see 
\cite{Bar93}) 
and can be modified here for our requirements.

Let $(M^{n,1},g)$ be a Lorentzian manifold. 
We consider the cone $\hat{M}$ of signature $(2,n-1)$ on $M^{n,1}$,
which is defined as 
\[\hat{M}:=(M\times\R_+\, ,\, \hat{g}:= r^2g-dr^2)\ .\]
The vector $r\cdot \partial_r$ is called the Euler vector of $\hat{M}$.
The $1$-level $M\times\{1\}$ of the cone $\hat{M}$ is naturally
isometric to the base manifold $M^{n,1}$ itself. 
We denote by 
$\tilde{X}$ the pullback of an arbitrary base vector field
$X\in\Gamma(M)$ to 
$\hat{M}$ through the projection $\pi$. Then we have the 
following rules for the Levi-Civita connection $\hat{\nabla}$ on the cone
\[\begin{array}{l}\hat{\nabla}_{\partial_r}\partial_r=0,
\qquad\hat{\nabla}_{\partial_r}\tilde{X}=\hat{\nabla}_{\tilde{X}}\partial_r
=\frac{1}{r}\tilde{X}\ ,\\[2mm]
\hat{\nabla}_{\tilde{X}}\tilde{Y}=\nabla_XY-
rg(X,Y)\partial_r\ .\end{array}\]

In case that $M^{n,1}$ is a spin manifold the cone
$\hat{M}$
is a spin manifold, too. Then we denote the spinor bundle of the
cone with $\hat{S}$. For $n$ even the restriction of
$\hat{S}$ 
to the $1$-level $M\times\{1\}$ of the cone is naturally isomorphic to the 
spinor bundle $S$ on the base manifold $M^{n,1}$ by a map
\[\Phi\ :\ S\ \cong\ \hat{S}|_{M\times\{1\}}
\]
with $\Phi(X\cdot\varphi)=X\cdot\Phi(\varphi)$ for all 
$X\in TM^{n,1}$.
Similar, if $n$ is odd, there are isomorphisms 
$\Phi_{\pm}:S\cong\hat{S}^{\pm}|_{M\times\{1\}}$
for the restricted half spinor bundles
such that
\[-iX\cdot\Phi_+(\varphi)=\Phi_-(X\cdot\varphi)\ . \]
for all tangent vectors $X\in TM^{n,1}$.  
With respect to the metric $\hat{g}$ the
projection $\pi$
gives rise to a pullback
$\pi^*:\Gamma(\hat{S}|_{M\times\{1\}})
\to\Gamma(\hat{S})$ of spinor fields on the $1$-level to the cone.
Eventually, we denote by
$\Cal{K}_\lambda(M)$ the space of Killing spinors on $(M^{n,1},g)$ to the 
Killing number $\lambda$. 

\begin{theorem} \label{t31} (cf. \cite{Bar93} and \cite{Sem02}) Let
$(M^{n,1},g)$ be a 
Lorentzian manifold and $\hat{M}$ its cone with signature $(2,n-1)$. 
The following correspondences exist.
\begin{enumerate}   
\item 
The special Killing $p$-forms on $M^{n,1}$ to the
positive constant $c=p+1$ 
are in 1-to-1-correspondence with the parallel $(p+1)$-forms on the
cone $\hat{M}$. The correspondence is given by
\[\alpha\in\Omega^p(M)\qquad\mapsto\qquad r^pdr\wedge
\alpha-\frac{r^{p+1}}{p+1}d\alpha\in\Omega^{p+1}(\hat{M})\ .\]
\item If $M^{n,1}$ is spin and $scal=-n(n-1)$ then there are natural
isomorphisms
\[\begin{array}{ccc}
\Cal{K}_{\frac{i}{2}}(M)\oplus\Cal{K}_{-\frac{i}{2}}(M)&\cong&\
\Cal{K}_0(\hat{M})\\[1mm]\varphi&\mapsto&\
\hat{\varphi}:=\pi^*\circ\Phi(\varphi)\end{array}\qquad\quad\qquad\ 
\mbox{for}\ n\ \mbox{even\ \ and}\] 
\[\begin{array}{ccc}
\Cal{K}_{\pm\frac{i}{2}}(M)&\qquad\cong&\ 
\ \Cal{K}^{\pm}_{0}(\hat{M})\\[1mm]
\varphi&\qquad\mapsto&\quad \hat{\varphi}:=\pi^*\circ\Phi_\pm(\varphi)\
\end{array}\qquad\qquad\qquad \mbox{for}\ n\ \mbox{odd},\]
where $\Cal{K}^{\pm}_{0}(\hat{M})$
is the space of parallel $\pm$-half spinors on the cone. 
\end{enumerate}
\end{theorem} 

The Riemannian version of Theorem \ref{t31} is classical for the
application
to the case of real Killing spinors. The 
result for
Killing $p$-forms on Riemannian manifolds was established in
\cite{Sem02}. The proof for the correspondence here in case of imaginary 
Killing 
spinors $\varphi$ 
in Lorentzian geomtery is 
based on the observation that $\varphi$ is parallel with respect to the 
modified spinor connection $\tilde{\nabla}_\lambda$
coming from an affine connection, which takes
values in 
\[ i\R^{1,n-1}\oplus\frak{spin}(1,n-1)
\cong\frak{spin}(2,n-1)\subset
\it{Cliff}^{\C}_{1,n-1}\ .\]
We remark for the application of Theorem \ref{t31} that the metric $g$ 
on $M^{n,1}$ can be rescaled by a 
positive constant such that the positive constant $c$ to an
arbitrary special Killing 
$p$-form equals $p+1$ and the Killing number $\lambda$ 
to an arbitrary imaginary Killing spinor satisfies
$\lambda^2=-\frac{1}{4}$.\\ 

The spinor bundle $S^{n,2}$ on a time-oriented
pseudo-Riemannian spin
manifold
$(N^{n,2},h)$ of signature $(2,n-2)$ is equipped with an invariant
inner product $\langle\cdot,\cdot\rangle_{2,n-2}$ (cf. \cite{Bau81}). 
Similar
to the induced Dirac current of a spinor in Lorentzian
geometry, a spinor $\gamma\in\Gamma(S^{n,2})$ induces
a $2$-form $\alpha^2_{\gamma}$ on $N^{n,2}$
by the rule
\[h(\alpha^2_\gamma,X^2):=-i\langle
X^2\cdot\gamma,\gamma\rangle_{2,n-2}\qquad\mbox{for\ all}\ 
2\mbox{-forms}\
X^2.\]

In case that $\hat{M}$ is the cone over a Lorenztian spin manifold
$M^{n,1}$ the inner product $\langle\cdot,\cdot\rangle_{2,n-1}$
admits the property
\[\begin{array}{lccl}\langle\varphi,\psi\rangle&=&
-\langle\partial_r\cdot\Phi_-(\varphi),\Phi_+(\psi)\rangle_{2,n-1}
&\qquad\mbox{for}\ n\ \mbox{odd\quad and}\\[3mm]
\langle\varphi,\psi\rangle&=&
i\langle\partial_r\cdot\Phi(\varphi),\Phi(\psi)\rangle_{2,n-1} 
&\qquad
\mbox{for}\ n\ \mbox{even}\end{array}\]
on the $1$-level of $\hat{M}$,
where $\varphi$, $\psi$ are spinor fields on $M^{n,1}$.
Then the following relation is true.

\begin{lemma} \label{l32}
Let $\varphi\in\Gamma(S)$ be a spinor with Dirac current $V_\varphi$ 
on a Lorentzian spin
manifold $M^{n,1}$ and let 
$\hat{\varphi}$ be the corresponding
($\pm$-half) spinor
with 
associated $2$-form $\alpha^2_{\hat{\varphi}}$ on the 
cone $\hat{M}$. It holds 
$V_{\varphi}^\flat=\partial_r\ecke\alpha^2_{\hat{\varphi}}$  
on the $1$-level $M^{n,1}\subset\hat{M}$. \end{lemma}

{\sc Proof}: With respect to an orthonormal basis $e=(e_0,e_1,\ldots,e_n)$
with
$e_0=\partial_r$ in an arbitrary point of the $1$-level it  holds
\begin{eqnarray*}\partial_r\ecke\alpha^2_{\hat{\varphi}}
&=&-i\sum_{i<j}\langle e_ie_j\cdot\hat{\varphi},
\hat{\varphi}\rangle_{2,n-1}\cdot\partial_r\ecke e_i^*\wedge e_j^*\\
&=&-i\sum_{j=1}^n\langle\partial_re_j\cdot\hat{\varphi},\hat{\varphi}
\rangle_{2,n-1} e_j^*=-\sum_{j=1}^n\langle e_j\varphi,\varphi\rangle e_j^*
=V^\flat_\varphi\ .\end{eqnarray*}
\hfill$\Box$\\

The lemma also shows that $\alpha^2_{\hat{\varphi}}$
is non-trivial for all (half) spinors $\hat{\varphi}\neq 0$ on the cone
$\hat{M}$,
since the corresponding Dirac current 
$V_{\varphi}^\flat$ on $M^{n,1}$ is non-trivial. 

\section{Normal forms for skew-adjoint                               %
         operators in signature $(2,n-2)$}                           %
 
In this section we present a complete list of {\it normal forms} for 
skew-adjoint endomorphisms acting on the pseudo-Euclidean space
$\R^{2,n-2}$ of dimension $n$ and signature $(2,n-2)$. 
This list was established in \cite{Bou00}. 
Parallel $2$-forms on the cone of signature 
$(2,n-1)$ over a Lorentzian manifold correspond to parallel skew-adjoint
operators and are therefore distinguished
by normal forms of the list. This observation
will be the crucial point
in our 
description of Lorenztian geometries admitting imaginary Killing spinors
in the last section.

\begin{theorem} \label{t41} (cf. \cite{Bou00})
Let $\beta$ be an arbitrary $2$-form on
the pseudo-Euclidean space $\R^{2,n-2}$.
Then there exist vector spaces $V_i$ such that
$\R^{2,n-2}=\oplus_i V_i$ is an orthogonal direct sum and
the skew-adjoint endomorphism $b$, which corresponds to
$\beta$, satisfies $b(V_i)\subset V_i$ for all $i$. Moreover,
there is a basis $(e_{i_1},\ldots,e_{i_{r(i)}})$ for every $V_i$ such  
that the corresponding matrix for
the inner product and for $b$ is one pair of blocks as
it occurs in the lines of Table 1 below.
\end{theorem}

The basis of $\R^{2,n-2}$, in which
a skew-adjoint operator takes a normal form, is called
an {\it adapted basis}. There is always an orthogonal decomposition
$\R^{2,n-2}=E\oplus P$ to a skew-adjoint operator $b$
such that $E$ is Euclidean and $b$ preserves the decomposition.
We call the normal form to $b$ on $E$ an {\it Euclidean block}
and the normal form to $b$ on $P$ a {\it pseudo-Euclidean block}.

\begin{example} \begin{enumerate}
\item[a)] Let $\omega_o:=\sum_{i=1}^m
e^*_{2i-1}\wedge
e_{2i}^*$ be the standard (pseudo)-K{\"a}hler
form on $\R^{2,n-2}$, where $(e_1,\ldots,e_{2m})$ is the standard basis.
The normal form of the  skew-adjoint operator
corresponding to a multiple $\omega=\nu\cdot\omega_o$ of the K{\"a}hler
form
with respect to the adapted  
basis $(e_1,\ldots,e_{2m})$
is built up by one block of the form
$B_{II}(\nu)$ (pseudo-Euclidean block) and $(m-1)$ blocks of
the form $B(\nu)$
(Euclidean block).
\item[b)] A $2$-form $\omega=l_1^\flat\wedge l_2^\flat$
on $\R^{2,n-2}$,
where $l_1$ and $l_2$ are lightlike vectors, which span a totally
lightlike
plane,
corresponds as
skew-adjoint
operator
with
respect to
some adapted basis to a composition
of a
pseudo-Euclidean block of the form $B_{Ia}$ and
an Euclidean $0$-block of length $n-4$.
\item[c)] A $2$-form $\omega=l_1^\flat\wedge t_1^\flat$
on $\R^{2,n-2}$,
where $l_1$ is lightlike, $t_1$ is timelike and both vectors
are orthogonal, corresponds
as skew-adjoint operator with respect to
some adapted basis to a composition of a block $B_{Ib}$ and
a $0$-block of length $n-3$.
\end{enumerate}
\end{example}

Let $\hat{\varphi}$ be a spinor on the pseudo-Euclidean space
$\R^{2,n-1}$. There corresponds a $2$-form $\alpha^2_{\hat{\varphi}}$
to $\hat{\varphi}$ on $\R^{2,n-1}$ defined by the rule
\[(\alpha^2_\varphi,x^2):=-i\langle
x^2\cdot\varphi,\varphi\rangle_{2,n-2}\qquad\mbox{for\ all}\
x^2\in\Lambda^2(\R^{{2,n-1}^*}),\]
where $(\cdot,\cdot)$ denotes the induced inner product on
$\Lambda^2(\R^{{2,n-1}^*})$  
(cf. section 3).
The following statement is a version of Lemma \ref{l32}
considered in a single point only and the proof for it works the same
as there before.
\begin{table}\scriptsize
\[\begin{array}{lccc}\!\!\!\!\!\!\!\!\!\!\!\!\AS\mbox{signature}\
(p,q)&\ \quad\AS A=\mbox{
inner\ product}\quad \ &\qquad\AS B=\mbox{skew-adjoint\
operator}\qquad\\ \\ 
(0,1):&(1)&(0)&\\ \\
(0,2):&\left(\begin{array}{cc}
1&0\\
0&1
\end{array} \right)&
B(\nu)=\left( \begin{array}{cc}
0&-\nu\\     
\nu&0    
\end{array} \right)\qquad\quad &\nu\neq 0\\ \\ 
(1,0):&(-1)&(0)&\\ \\
(1,2):&
\left( \begin{array}{ccc}
0&0&-1\\
0&1&0\\
-1&0&0
\end{array} \right)
&\left( \begin{array}{ccc}
0&1&0\\
0&0&1\\    
0&0&0
\end{array} \right)&\\ \\
(1,1):&
\left( \begin{array}{cc}
0&1\\
1&0
\end{array} \right)
&\left( \begin{array}{cc}
\lambda&0\\
0&-\lambda
\end{array} \right)&\lambda\neq 0\\ \\
(2,2):&\left( \begin{array}{cccc}
0&0&0&-1\\
0&0&1&0\\ 
0&1&0&0\\
-1&0&0&0
\end{array} \right)
&
B_{Ia}=\left( \begin{array}{cccc}
0&0&1&0\\
0&0&0&1\\
0&0&0&0\\
0&0&0&0
\end{array} \right)\\ \\
(2,1):&\left( \begin{array}{ccc}
0&0&1\\
0&-1&0\\ 
1&0&0
\end{array} \right)
&
B_{Ib}=\left( \begin{array}{ccc}
0&1&0\\
0&0&1\\
0&0&0
\end{array} \right)\quad\quad \\ \\
(2,4):&\left( \begin{array}{ccc}
0&0&-I_2\\
0&I_2&0\\ 
-I_2&0&0
\end{array} \right)
&\left( \begin{array}{ccc}
0&I_2&0\\
0&0&I_2\\
0&0&0
\end{array} \right)\\ \\
(2,0):&\left( \begin{array}{cc}
-1&0\\
0&-1
\end{array} \right)&
B_{II}(\nu)=\left( \begin{array}{cc}
0&-\nu\\
\nu&0
\end{array} \right)\quad\quad &\nu\neq 0\\ \\
(2,2):&\pm\left( \begin{array}{cccc}
0&0&0&-1\\
0&0&1&0\\
0&1&0&0\\
-1&0&0&0
\end{array} \right)
&\ B^{\pm}_{IIa}\!=\!\!\left(\!\!\!\! \begin{array}{cccc}
0&\!\!\!\! -\nu_\pm&\!\!\!\! 1&\!\!\!\! 0\\
\nu_\pm&\!\!\!\! 0&\!\!\!\! 0&\!\!\!\! 1\\
0&\!\!\!\! 0&\!\!\!\! 0&\!\!\!\! -\nu_\pm\\
0&\!\!\!\! 0&\!\!\!\! \nu_\pm&\!\!\!\! 0
\end{array}\!\!\!\!\! \right)&\nu_\pm\neq 0\\ \\
(2,4):&\left( \begin{array}{ccc} 
0&0&-I_2\\
0&I_2&0\\
-I_2&0&0 
\end{array} \right)
&\quad  \left(\!\!\!\! \begin{array}{ccc}
\scriptsize\begin{array}{cc}
0&-\nu\\
\nu&0
\end{array}&\!\!\!\!\! I_2&\!\!\!\!\!0\\
0&\!\!\!\!\!\scriptsize\begin{array}{cc}
0&-\nu\\
\nu&0
\end{array}&\!\!\!\!\! I_2\\
0&\!\!\!\!\! 0&\!\!\!\!\!\scriptsize\begin{array}{cc}
0&-\nu\\
\nu&0
\end{array}
\end{array}\!\!\!\!\!\! \right)\quad &\nu\neq 0\\ \\
(2,2):&\left( \begin{array}{cc}
0&I_2\\  
I_2&0
\end{array} \right)
&\left( \begin{array}{cc}
\lambda I_2&0\\ 
0&-\lambda I_2
\end{array} \right)&\lambda\neq 0\\ \\
(2,2):&\qquad\left( \begin{array}{cccc}
0&0&0&-1\\
0&0&1&0\\
0&1&0&0\\
-1&0&0&0
\end{array} \right)\qquad&\qquad
\left( \begin{array}{cccc}
\lambda&0&1&0\\
0&-\lambda&0&1\\
0&0&\lambda&0\\
0&0&0&-\lambda
\end{array} \right)\qquad&\lambda\neq 0\\ \\
(2,2):&\qquad\left( \begin{array}{cccc}
0&0&1&0\\
0&0&0&-1\\
1&0&0&0\\
0&-1&0&0
\end{array} \right)\qquad&\quad\
B_{IIb}=\!\left(\!\! \begin{array}{cccc}
\xi\!&-\nu\!&0\!&0\\
\nu\!&\xi\!&0\!&0\\
0\!&0\!&-\xi\!&\nu\\
0\!&0\!&-\nu\!&-\xi
\end{array} \!\!\right)\ \ &\xi,\nu\neq 0\\[5mm]
\end{array}\]   
\caption{{\footnotesize These are the building
blocks for
the normal forms of skew-adjoint operators in signature $(2,n-2)$.
The matrices in the first column (denoted by $A$) indicate
an inner product (of index $s\leq 2$) with respect to some
basis and the matrices in the second
column (denoted by $B$) are skew-adjoint endomorphisms with respect to the
inner product in column $A$ and the chosen basis.}}
\end{table}

\begin{lemma} \label{l43} Let $\hat{\varphi}$ be a
($\pm$-half) spinor on $\R^{2,n-1}$ and $T\in\R^{2,n-1}$
an arbitrary unit timelike vector.
The $1$-form
$\alpha_{T,\hat{\varphi}}:= T\ecke\alpha^2_{\hat{\varphi}}$
is dual to the associated vector
induced by the spinor $\varphi$ on the Minkowski space
$T^{\bot}\subset\R^{2,n-1}$, which corresponds naturally
to $\hat{\varphi}$.
\end{lemma} 


The lemma imposes a condition on the nature 
of a $2$-form induced by a spinor in signature $(2,n-2)$, since
the associated vector to a non-trivial spinor on the Minkowski space 
is not arbitrary, but causal. 
With some simple calculations we
can sort out the normal
forms for skew-adjoint operators corresponding to $2$-forms,
which do not satisfy the condition imposed by Lemma \ref{l43}, 
and therefore can not be induced by a spinor.
 
\begin{co} \label{c44} Let $\omega$ be a $2$-form in signature
$(2,n-2)$
such that the covector $T\ecke\omega$ is causal for every
timelike vector $T\in\R^{2,n-2}$. 
\begin{enumerate}
\item If there is a timelike $T$ such that
$T\ecke\omega$ is lightlike then the normal form corresponding to 
$\omega$ is a   
composition of a pseudo-Euclidean block of the form
$B_{Ia}$ or $B_{Ib}$ with an Euclidean $0$-block. 
\item If $T\ecke\omega$ is timelike for all timelike $T$ then
the normal form of $\omega$ 
is a composition
of $B_{II}$, $B^{\pm}_{IIa}$ ($\nu_+>0$,
$\nu_-<0$) or
$B_{IIb}$ ($\nu^2\geq\xi^2$) with an
Euclidean  block consisting of blocks of the form $B(\nu)$
and/or a $0$-block. 
\end{enumerate}
\end{co}

With {\it stabilizer} of a skew-adjoint operator (resp. $2$-form)
we mean in the following the subgroup of the 
(pseudo)-orthogonal group, which leaves the operator (resp. the $2$-form)
invariant under conjugated action. A simple consideration shows the
following fact.

\begin{lemma} \label{l45}
The stabilizer of a normal form, which is built from
a pseudo-Euclidean block of the form $B_{IIa}^\pm(\nu_\pm)$ or
$B_{IIb}(\nu,\xi)$
and some Euclidean block of length $n-4$ is included
in $U(1,1)\times \SO(n-4)$ for all eigenvalues $\nu_\pm,\nu$ and
$\xi\neq 0$.
\end{lemma}

\begin{definition} Let $\R^{2,n-2}$ be the pseudo-Euclidean
space of signature $(2,n-2)$ and $\omega\in\Lambda^2\R^{{2,n-2}^*}$ be a
non-trivial
2-form. We say that $\omega$ is of 
\begin{enumerate}
\item[$\bullet$]
Type $(I_a)$\quad if\,  $\omega=l_1^\flat\wedge l_2^\flat$ for some
vectors $l_1$ and $l_2$, which span a totally lightlike plane. 
\item[$\bullet$] Type $(I_b)$\quad if\,  $\omega=l_1^\flat\wedge
t_1^\flat$
for some
lightlike vector $l_1$ and a $l_1$-orthogonal timelike vector $t_1$.
\item[$\bullet$] Type $(II_a)$ or 
K{\"a}hler Type\quad if\, $\omega$ is a non-trivial multiple of the
K{\"a}hler form.
\item[$\bullet$]
Type $(II_b)$\quad if\,  there exists a non-trivial  Euclidean subspace
$E$
in
$\R^{2,n-2}$
such that $\omega$ restricted to $E$ vanishes and $\omega$ is
a (pseudo)-K{\"a}hler form on the orthogonal complement of $E$ in
$\R^{2,n-2}$.
\end{enumerate}
\end{definition} 

Lemma \ref{l45} makes clear, which stabilizers of the normal forms
occuring in Corollary \ref{c44} are maximal.

\begin{co} \label{c47}
A $2$-form $\omega$ on $\R^{2,n-2}$, which is of Type $(I_a)$,
$(I_b)$ , $(II_a)$
or  $(II_b)$,
is exclusively distinguished by its properties that
\begin{enumerate}
\item
the covector $T\ecke\omega$ is causal for every timelike vector 
$T\in\R^{2,n-2}$ and 
\item
its stabilizer $S_\omega$ in $\SO(2,n-2)$
is maximal, in the sence that there is no non-trivial $2$-form
satisfying the first property, whose stabilizer properly contains
$S_\omega$.
\end{enumerate}
\end{co}

We observe
that the stabilizer of a $2$-form of Type $(I_a)$ and Type $(I_b)$ acts
indecomposable but reducible on $\R^{2,n-2}$, i.e. there exist
non-trivial and invariant subspaces of $\R^{2,n-2}$, but the inner
product is degenerate on all of them. The
stabilizer of a K{\"a}hler type form is $U(1,m-1)$ and acts irreducible
on $\R^{2,2m-2}$. The stabilizer of a form of Type $(II_b)$ acts
decomposable on $\R^{2,n-2}$.

\section{Imaginary Killing spinors}  
With the construction of the cone in section 3 and the normal
form classification for skew-adjoint operators in signature $(2,n-2)$
(coming from a spinor) in the last section we are now in the position to
discuss
a geometric description of Lorentzian manifolds admitting imaginary
Killing spinors.\\

In the following, the metric $g$ on the Lorenztian spin manifold
$M^{n,1}$ will be scaled such that the
Killing number $\lambda$ to any Killing spinor 
satisfies $\lambda^2=-\frac{1}{4}$.
We start with a proposition, which characterizes the cone $\hat{M}$
of a Lorentzian manifold $M^{n,1}$ with imaginary Killing spinor and  
indicates the different cases that are to be considered
for the geometry of the base manifold $M^{n,1}$. We remark 
that the normal form corresponding to a parallel $2$-form on the cone
is in every point the same.

\begin{prop} \label{p51} Let $(M^{n,1},g)$ be a Lorentzian spin manifold
admitting
an imaginary Killing spinor. Either there exists
a parallel $2$-form $\omega$ of Type $(II_b)$ on the cone $\hat{M}$ or 
there exists at least one parallel (half) spinor $\hat{\varphi}$ on
$\hat{M}$
such that the induced 
parallel $2$-form $\omega=\alpha^2_{\hat{\varphi}}$ is of
Type $(I_a)$, $(I_b)$ or $(II_a)$.
\end{prop}

{\sc Proof}: Let $\psi$ be an imaginary Killing spinor
on $M^{n,1}$. According to
Corollary
\ref{c44} the normal form of the 
skew-adjoint endomorphism corresponding to $\alpha^2_{\hat{\psi}}$ is a
composition with one block of the form $B_{Ia}$, $B_{Ib}$,
$B_{II}$, $B_{IIa}^\pm$ or $B_{IIb}$. In case that the normal form is
built with a block of the form $B_{Ia}$ or $B_{Ib}$ the
parallel $2$-form $\omega=\alpha^2_{\hat{\psi}}$
is of Type $(I_a)$ or $(I_b)$.\\

In the other cases there exists a biggest number $s>0$ such
that the stabilizer of the
normal form to $\alpha^2_{\hat{\psi}}$ is included in 
$U(1,s-1)\times\SO(n-2s)$. This group
includes the holonomy group of the cone $\hat{M}$.
In case that $2s=n$ the cone $\hat{M}$ is a 
K{\"a}hler spin manifold. Moreover, since 
$V^\flat_\psi=\partial_r\ecke\alpha^2_{\hat{\psi}}$ is everywhere timelike
(Corollary \ref{c44}), the base $M^{n,1}$ is Einstein and hence the cone
is Ricci-flat. This implies that
there exists a parallel (half) spinor $\hat{\varphi}$, which
induces a K{\"a}hler form
on the cone. If $2s<n$ there exists a parallel $2$-form 
$\omega$ of Type $(II_b)$.\hfill$\Box$\\
   
We discuss now a description of the Lorentzian geometries on the base
manifold $M^{n,1}$ with imaginary Killing spinor 
according to the cases $(I_a)$, $(I_b)$ and $(II_a)$ that occur in
Proposition \ref{p51}. 

\subsection{Type $(I_a)$} 
In this case there exists a parallel (half) spinor $\hat{\varphi}$
on the cone $\hat{M}$, which induces a parallel $2$-form
$\omega\neq 0$ that is locally of the form $l_1^\flat\wedge
l_2^\flat$
for some lightlike vector fields $l_1$ and $l_2$, which span a totally
lightlike plane. 
The dual $V_\varphi^\flat$ of
the Dirac current of the imaginary Killing spinor $\varphi$, which
corresponds to $\hat{\varphi}$
on $\hat{M}$, is equal to
$\partial_r\ecke\omega$,
which shows that the Dirac current $V_\varphi$ is everywhere lightlike.\\

There is a known description of Lorentzian metrics
admitting {\it twistor spinors} with lightlike Dirac current. 
We call a Lorentzian space admitting a lightlike parallel
vector field a {\it Brinkmann space}.
Two spinor fields on (pseudo)-Riemannian spaces are said to be 
{\it conformally equivalent} 
if there exists a conformal diffeomorphism, which identifies both spinor
fields. In particuar, it holds

\begin{prop} \label{p52} (see \cite{BL02}) Let $\varphi$ be a spinor
field,
which satisfies the twistor equation 
$\nabla^S_X\varphi+\frac{1}{n}X\cdot D\varphi=0$ for all vector fields
$X$, such that the Dirac current $V_\varphi$ is a lightlike
Killing vector field on $M^{n,1}$. If
$Ric(V_\varphi,V_\varphi)=0$ then $\varphi$ is locally conformally
equivalent to a parallel spinor on a Brinkmann space.
\end{prop}

This gives rise to

\begin{prop} 
Let $\varphi$ be an imaginary Killing spinor on $M^{n,1}$ such that
$\alpha^2_{\hat{\varphi}}$ on $\hat{M}$ is of Type $(I_a)$.
Then $\varphi$ is locally conformally equivalent to a parallel spinor
on a Brinkmann space.
\end{prop}

{\sc Proof}: From Proposition \ref{p21} and Lemma \ref{l21} we know that
$Ric(V_\varphi)=\rho V_\varphi$ for some real function $\rho$. 
Then we can apply Proposition \ref{p52} to prove the result.
\hfill$\Box$


\subsection{Type $(I_b)$} 
There exists a parallel (half) spinor $\hat{\varphi}$
on the cone $\hat{M}$, which induces a parallel $2$-form $\omega$ of Type
$(I_b)$. In this situation it holds

\begin{lemma} \label{l51} The function
$f_\varphi:=\sqrt{-g(V_\varphi,V_\varphi)}$ 
to the imaginary Killing spinor $\varphi$ on $M^{n,1}$ satisfies 
\begin{enumerate}
\item $Hess(f_\varphi)=f_\varphi\cdot g$, i.e. $gradf_\varphi$ is a
conformal gradient field, and  $f_\varphi^2=g(gradf,gradf)$,
\item $gradf_\varphi\neq 0$ and
$f_\varphi\neq 0$ on
disjoint dense subspaces in $M^{n,1}$.
\end{enumerate}
\end{lemma}

{\sc Proof}: The $2$-form $\omega$ can be written as $rdr\wedge
V^\flat_\varphi-\frac{r^2}{2}dV_\varphi^\flat$ 
(Theorem \ref{t31}). 
The 2-dimensional parallel subbundle $E_{\omega}\subset
TM^{n,1}$, which corresponds to the indecomposable $\omega$ is
degenerate and there is a unique parallel lightlike direction in
$E_{\omega}$. In particular, there exists a parallel lightlike vector
field
$l_1$ on $\hat{M}$. Moreover, we can find locally a timelike vector $t_1$
of constant length such that $\omega=l_1^\flat\wedge t_1^\flat$. We 
choose the parallel lightlike field $l_1$ with the scaling 
$\hat{g}(t_1,t_1)=-1$.
Since $l_1$ is parallel, there is a 
unique function $f$ on $M^{n,1}$ such that 
$l_1^\flat=f dr-rdf$ and $f^2=g(gradf,gradf)$. The function $f$ is
a special Killing $0$-form on $M^{n,1}$, i.e. $gradf$ is a conformal
gradient field. Since neither $dr$ nor the lift of $df$ to the
cone are parallel, it is $f\neq 0$ and $gradf\neq 0$ on a
dense subset of $M^{n,1}$.\\

We calculate the function $f$ with respect to $\varphi$.
The local field $t_1^\flat$ is given by $t_1^\flat=Adr+u$,
where $A$ is a function and $u$ a $1$-form on $M^{n,1}$.
It follows that $V_\varphi^\flat=fu+Adf$ and 
\[g(V_\varphi,V_\varphi)=f^2g(u,u)+A^2g(df,df)+2fA\cdot g(u,df).\]
Since $g(u,u)=-1+A^2$ and $g(u,df)=-Af$, we can conclude
$f^2=-g(V_\varphi,V_\varphi)$, which shows that $f_\varphi$
has the claimed properties.\hfill$\Box$\\

The assertions of Lemma \ref{l51} imply together with
Proposition \ref{p22} that $M^{n,1}$ is Einstein and 
$gradf_\varphi$ is a non-homothetic conformal gradient field.
There is a known description of (pseudo)-Riemannian Einstein metrics
admitting such conformal fields. In particluar, there is

\begin{prop} \label{p54} (cf. \cite{KR97}) Let $(M^{n,1},g)$ be a
Lorentzian
Einstein
space
admitting
a non-constant solution $f$ of the equation $Hess(f)=l\cdot g$ for some
function $l$. Then, 
in a neighborhood of any point with $v:=g(grad(f),grad(f))\neq 0$,
the metric $g$ is a warped product 
$\varepsilon\cdot dt^2+f^{'2}(t)k$, where
$\varepsilon:=sign(v)$, $k$
is an Einstein
metric and $f$ satisfies $f^{''2}+\frac{\varepsilon scal_g}{n(n-1)}
f^{'2}=\frac{\varepsilon
scal_k}{(n-1)(n-2)}$.
\end{prop}

This leads to

\begin{prop}
Let $\varphi$ be an imaginary Killing spinor on $M^{n,1}$ such that
$\alpha^2_{\hat{\varphi}}$ on $\hat{M}$ is of Type $(I_b)$.
Then, in a neighborhood of any point with $V_\varphi$ timelike,
the metric $g$ is a warped product of the form $dt^2+f^{2}k$, where
$k$ is a Lorentzian Einstein metric admitting a Killing spinor to
the Killing number 
\begin{enumerate} 
\item $\lambda_k=0$ and $f=\exp t$
\item $\lambda_k=\frac{1}{2}$ and $f=\sinh t$ or
\item $\lambda_k=\frac{i}{2}$ and $f=\cosh t$.
\end{enumerate} 
\end{prop}

{\sc Proof}: The function $f_\varphi=\sqrt{-g(V_\varphi,V_\varphi)}$
satisfies the assumptions
of Proposition \ref{p54}. Since $f_\varphi^2>0$, the warping function 
$f=f^{'}_\varphi$
must
solve 
the ordinary differential equation $f^{'2}-n(n-1)f^2=scal_k$. 
There are three different solutions $f=\exp t$, $\cosh t$ and
$\sinh t$ according to the values $scal_k=0,\pm (n-1)(n-2)$.
In each case the imaginary Killing spinor $\varphi$ induces a
Killing spinor to the Killing number $\frac{scal_k}{(n-1)(n-2)}$ on 
the space with Einstein metric $k$ (cf. \cite{Boh98}).\hfill$\Box$


\subsection{Type $(II_a)$}
\begin{lemma} \label{l52}
Let $(M^{n,1},g)$ be a Lorenztian Einstein manifold with 
a Killing vector $V$ such that $g(V,V)=-1$ is constant and $V\ecke
\Cal{W}=0$ ($\Cal{W}$ Weyl tensor). Then
the operator $J$ defined by $J(X):=\nabla_XV$ on $TM$ satisfies
\begin{enumerate}
\item
$J(V)=0$\quad and\quad $J^2(X)=\frac{scal}{n(n-1)}\big(X+g(V,X)V\big)$
\item
$(\nabla_XJ)(Y)=\frac{scal}{n(n-1)}\big(g(V,Y)X-g(X,Y)V\big)$.
\end{enumerate}
\end{lemma}   

{\sc Proof}: Because $V$ is Killing with constant length, it follows
$\nabla_VV=0$
and $g(\nabla_XV,\nabla_YV)=\Cal{R}(V,X,Y,V)$, where $\Cal{R}$
denotes the Riemannian curvature tensor. It is
$\Cal{R}=\Cal{W}+g\star L$, where $L=\frac{1}{n-2}\big(\frac{scal}{2(n-1)}
g-Ric\big)$ is the Schouten tensor and $\star$ denotes the Kulkarni-Nomizu
product (cf. \cite{Bes87}).
Then from $V\ecke\Cal{W}=0$ we obtain
\begin{eqnarray*}
g(J^2(X),Y)&=&-g(J(X),J(Y))\\&=&-g(V,Y)L(X,V)-g(V,X)L(Y,V)-L(X,Y)
+g(X,Y)L(V,V)\ .\end{eqnarray*}  The relation for $J^2$ follows
immediately, since for
$M^{n,1}$ an Einstein space it holds $L=-\frac{scal}{2(n-1)n}g$.
Moreover, it is $g(\nabla_{e_k}\nabla_{e_i}V,e_j)=
\Cal{R}(e_i,e_j,e_k,V)$ for all $i,j,k\in\{1,\ldots,n\}$ in 
$p\in M^{n,1}$ arbitrary, 
where $(e_1,\ldots,e_n)$ is a local parallel frame in $p$.
Then 
\begin{eqnarray*}g((\nabla_{e_i}J)(e_k),e_l)&\ 
=\quad&\Cal{R}(e_k,e_l,e_i,V)\\
&\ =\quad&g(e_k,e_i)L(e_l,V)+g(e_l,V)L(e_k,e_i)\\
&\ \quad - &g(e_k,V)L(e_l,e_i)
-g(e_l,e_i)L(e_k,V)\ ,\end{eqnarray*}
which showes the identity for $\nabla J$ in an arbitrary
point $p$ of
$M^{n,1}$.
\hfill$\Box$\newpage
 
A Lorentzian manifold $(M^{n,1},g,V)$ with $V$ a timelike Killing vector
of constant length such that
the operator $J=\nabla V$ 
satisfies the both properties of Lemma \ref{l52}, is called a Lorentzian
Sasaki manifold. It is well-known that a Sasaki
structure $(V,J)$ on $M^{n,1}$
corresponds to a K{\"a}hler structure on the cone $\hat{M}$
(cf. \cite{Bar93} and \cite{Bau00}).

\begin{prop} 
A Lorenztian spin manifold $(M^{n,1},g)$ with an imaginary Killing spinor
$\varphi$, whose Dirac current $V_\varphi$ is timelike 
and has constant
length, is a Lorentzian Einstein-Sasaki manifold. This is exactly
the case when the lift $\hat{\varphi}$ induces a K{\"a}hler form
on the cone $\hat{M}$. 
\end{prop}  
 
{\sc Proof}: We have only to show that $V_\varphi\ecke\Cal{W}=0$ 
on $M^{n,1}$ and then apply Lemma \ref{l52}. With the identity
$\Cal{W}(\eta)\cdot\varphi=0$ 
(Proposition \ref{p21})
and the relation
$X\cdot\eta=-X\ecke\eta+X^\flat\wedge\eta$ in the Clifford algebra,
where $X$ denotes a vector and $\eta$ a $2$-form,
we obtain
\[
\Cal{W}(V_\varphi,X,Y,Z)=\langle\varphi,\Cal{W}(X,Y,Z)\cdot\varphi
\rangle=\langle\varphi,Z^\flat\wedge\Cal{W}(X,Y)
\cdot\varphi\rangle\in\R\ \
\mbox{for\ all}\ X,Y,Z\in TM.\]
But $\langle\varphi,\rho^3\cdot\varphi\rangle\in i\R$ 
for all $3$-forms $\rho^3$, and therefore $V_\varphi\ecke\Cal{W}=0$.
\hfill$\Box$\\

We summarize the different cases to

\begin{theorem} \label{t53} Let $(M^{n,1},g)$ be a Lorentzian spin
manifold
with imaginary Killing spinor $\varphi$.
\begin{enumerate}
\item
If $M^{n,1}$ is not Einstein then $M^{n,1}$ 
is locally conformally equivalent 
to a Brinkmann space with parallel spinor.
\item
If $g(V_\varphi,V_\varphi)$ is constant then
\begin{enumerate}
\item[i)] $g(V_\varphi,V_\varphi)=0$ and $M^{n,1}$ is
locally conformally equivalent to a Brinkmann space with parallel spinor
or
\item[ii)] $g(V_\varphi,V_\varphi)<0$ and $M^{n,1}$
is a Lorentzian Einstein-Sasaki manifold.
\end{enumerate}
\item
If the cone $\hat{M}$ is indecomposable and $V_\varphi$ is timelike 
then $M^{n,1}$ is either
\begin{enumerate}
\item[i)] locally conformally equivalent to a Brinkmann space with
parallel spinor,
\item[ii)] locally a warped product of the form $dt^2+f^2k$,
where $k$ is a Lorentzian Einstein metric admitting a
Killing spinor and $f=\exp t$, $\cosh t$ or $\sinh t$ or
\item[iii)] a Lorentzian Einstein-Sasaki space (and the cone
$\hat{M}$ is irreducible).
\end{enumerate}
\item If $V_\varphi$ changes the causal type then the set
$Z_\varphi\subset
M^{n,1}$, where $V_\varphi$ is lightlike, is a hypersurface and 
$M^{n,1}\backslash Z_\varphi$ admits locally a warped product structure as
in {\it
3. ii)}.
\end{enumerate}
In case that the metric $g$
does not belong to one of those listed in {\it{3.}} then either
$V_\varphi$
changes
the causal type or there is a parallel $2$-form of Type $(II_b)$ on
the cone $\hat{M}$. 
\end{theorem} 

\begin{example}
Let $H^{n,1}:=\{ x\in\R^{2,n-1}\, :\, ||x||^2=-1\}\subset
\R^{2,n-1}$ be the pseudo-hyperbolic space of signature $(1,n-1)$ with
negative scalar curvature $scal=-n(n-1)$. The space $H^{n,1}$ is
geodesically
complete, time-orientable and spin. The cone over $H^{n,1}$ is
an open subset of $\R^{2,n-1}$. Each parallel (half) spinor on
$\R^{2,n-1}$
restricted to $H^{n,1}$ gives rise to an imaginary Killing spinor.
It is not difficult to see that every generic type $(I_a)$, $(I_b)$,
$(II_a)$ and $(II_b)$ is realized by a $2$-form, which comes from a
parallel (half) spinor on $\R^{2,n-1}$ and thus belongs to an imaginary 
Killing spinor on $H^{n,1}$. This means that there are examples
of imaginary Killing spinors on $H^{n,1}$, $n\geq 3$, for each 
case where the Dirac
current
is everywhere lightlike $(I_a)$, changes the causal type $(I_b)$, timelike
with constant length $(II_a)$ or everywhere timelike with non-constant
length.
\end{example}

\begin{remark}
\begin{enumerate}
\item
In fact, there exist examples of imaginary Killing spinors on
non-Einstein spaces, which are generated by an appropriate
conformal change of certain non-Einstein Brinkmann spaces 
with parallel spinors (cf. \cite{Boh98}).
\item
Partial structure results and examples for Lorenztian metrics with
parallel or real Killing spinors are known (cf. e.g. \cite{Bry00},
\cite{Boh98}).
The warped product structure in case of Type $(I_b)$ then provides 
a concrete
construction principle for imaginary Killing spinors on Lorenztian
spaces.

The complete result in \cite{KR97} for the description of Einstein
spaces with conformal gradient fields does not apply when the field
changes the causal type. As consequence, Theorem
\ref{t53}
does not describe Lorentzian metrics with imaginary Killing spinors
when the Dirac current changes the causal type.
\item
There is a construction principle for Lorenztian Einstein-Sasaki
spin spaces. They appear as $S^1$-fiber bundles over Riemannian 
K{\"a}hler-Einstein spin spaces of negative scalar curvature
(cf. \cite{Bau00}).
\item
In case that there exists a parallel $2$-form of Type $(II_b)$
the cone $\hat{M}$ is decomposable. Different from the Riemannian case,
this does not imply that the cone is flat, even if the base $M^{n,1}$
is geodesically complete. The geometry of the base
$M^{n,1}$
with imaginary Killing spinor in this case remains to be investigated
and is subject of a forthcoming paper
\end{enumerate}

\end{remark}

{\it Acknowledgment}: I would like to thank Thomas Neukirchner, Helga
Baum and Jose Figueroa-O'Farrill for many helpful discussions and
comments.

{\footnotesize


}
{\normalsize
{\sc
University of Edinburgh, School of Mathematics,\\ 
JCMB--King's Buildings, EH9 3JZ Edinburgh, Scotland\\}
{\it E-mail address:} {\tt felipe@maths.ed.ac.uk}}

\end{sloppypar}
\end{document}